\documentclass[a4paper,10pt]{elsarticle}
\usepackage{amsmath}
\usepackage{amssymb}
\usepackage{amsthm}
\usepackage{amsfonts}
\usepackage{graphicx}
\usepackage{rotating}
\usepackage{setspace}
\usepackage{url}
\usepackage{epsfig,tikz,multirow}
\usepackage{epstopdf}
\usepackage{subfig}
\usepackage{color}
\newtheorem{theorem}{Theorem}[section]

\newtheorem{lemma}{Lemma}[section]

\theoremstyle{definition}

\newtheorem{example}{Example}[section]

\numberwithin{equation}{section}

\journal{...}

\begin{document}

\begin{frontmatter}
\title{Positive Solutions for  Nonlinear Elliptic Equations  Depending on a Parameter with Dirichlet Boundary Conditions}

\author[add1]{Seshadev Padhi\corref{cor1}}
\author[add2]{John R Graef\corref{cor2}}
\ead{Your email please }
\author[add1]{Ankur Kanaujiya}
\address[add1]{Department of Mathematics, Birla Institute of Technology, Mesra, Ranchi - 835215, India.}
\address[add1]{Department of Mathematics, University of Tennessee at Chattanooga, Chattanooga, USA.}

\begin{abstract}
   We prove new results on the existence of positive radial solutions of the elliptic equation $-\Delta u= \lambda h(|x|,u)$   in an annular domain in $\mathbb{R}^{N}, N\geq 2$. Existence of positive radial solutions are determined under the conditions that the nonlinearity function $h(t,u)$ is either  superlinear or sublinear growth in $u$  or satisfies some upper and lower inequalities on $h$.  Our discussion is based on a fixed point theorem due to  a revised version of a fixed point theorem of Gustafson and Schmitt.
  \end{abstract}
\end{frontmatter}

\section{Introduction}

\ \ \ \ We consider the existence of positive solutions for the elliptic boundary value problem
 \begin{equation}\label{e1}
 \begin{cases}
  -\Delta v  = \lambda h(|x|,v) \,\, \text{in}\,\, \Omega  \\
 v  = 0  \,\, \text{on}\,\, \partial \Omega,
\end{cases}
  \end{equation}
where $\Omega = \{x \in \mathbb{R}^{N};r_{1}<|x|<r_{2}\}$ with $0<r_{1}<r_{2}, N \geq 2$, $h:[r_{1},r_{2}]\times \mathbb{R}^{+} \to \mathbb{R}^{+}$ is a continuous function, $h(t,0)=0$ and $\lambda>0$ is a real number. \\

\ \ \ \ In recent years, study on the existence of positive solutions for  elliptic equations of the form (\ref{e1}) and its various versions with Dirichlet and/or Neumann type boundary conditions have been given a serious attention. This is evident from the works in \cite{r1, r3,r4,r5,r6,r7,r8,r9,r10,r11, r14, r18,r19, r21,r22,r23,r24,r25,r26,r27,r29,r31,r32,r33}. \\

\ \ \ \ Elliptic equations of the form
 \begin{equation}\label{e2}
\begin{cases}
  -\Delta v  =\lambda h(v),  \,\, v  \in \Omega  \\
\ \ \ \ v  =0 \,\,  \text{on}\,\, \partial \Omega,
\end{cases}
  \end{equation}
 have been studied in \cite{r4,r7,r11,r24,r27,r28,r30}. In \cite{r30}, Shivaji proved some sharp conditions on the uniqueness of positive solutions of (\ref{e2}). Dancer and Schmitt \cite{r7} obtained some necessary conditions for the existence of positive solutions of (\ref{e2}), whose supremum norm bears a certain relationship to zeros of the nonlinearity $h$. Maya and Shivaji \cite{r28} used sub-super solution method to find the existence and nonexistence of positive solutions of (\ref{e2}). The results due to Maya and Shivaji \cite{r28} were extended by Perera \cite{r29} to  quasilinear elliptic problems using variational arguments. The authors in \cite{r7,r27,r28,r29,r30} considered a general bounded domain $\Omega$ in $\mathbb{R}^{N}, \, N \geq 2$. On the other hand, Lin \cite{r24} studied the existence and multiplicity of positive radial solutions of (\ref{e2}), where $\Omega$ is an annular domain of $\mathbb{R}^{N}\, N \geq 2$. Lin \cite{r24} proved that if $h(u)>0$ for $u \geq 0$ and $\lim_{u \to \infty}\frac{h(u)}{u}=\infty$, then there exists $\lambda^{*}>0$ such that there are at least two positive radial solutions for each $\lambda \in (0, \lambda^{*})$, at least one for $\lambda =\lambda^{*}$ and  none for $\lambda > \lambda^{*}$. If $h(0)=0$ and $\lim_{u \to 0}\frac{h(u)}{u}=1$ and $uh^{\prime}(u)>(1+\epsilon)h(u)$  for $u>0$ and $\epsilon>0$, then there exists a variational solution of (\ref{e2})  for $\lambda \in (0,\lambda_{1})$, where $\lambda_{1}$ is the least eigen value of $-\Delta$. Further, if $h(0)=0$, $\lim_{u \to 0}\frac{h(u)}{u}=0$ and $\lim_{u \to \infty}\frac{h(u)}{u}=\infty$, there exists at least one positive radial solution of (\ref{e2}) for any $\lambda>0$.\\

\ \ \ \ Erbe and Wang \cite{r9} used cone expansion and compression theorem to study the existence of positive solutions of  the elliptic equation
 \begin{equation}\label{e3}
\begin{cases}
  -\Delta u=\lambda g(|x|)f(u) \,\, u \in \Omega  \\
v=0 \,\, x \in \partial \Omega,
\end{cases}
  \end{equation}
where $\Omega \subset \mathbb{R}^{N}$ is an annular domain and $N \geq 1$. Similar equations have also been studied in \cite{r1,r8,r23,r33} using Mountain Pass theorem, Shooting method and different fixed point theorems.

\ \ \ \ Let $N=2$. If we set $r=r_{2}\left(\frac{r_{1}}{r_{2}}\right)^{t}$ and $u(t)=v(r)$, then (\ref{e1}) can be transformed to the boundary value problem (BVP in short)
 \begin{equation}\label{e4}
\begin{cases}
 u^{\prime\prime}(t)+\lambda q(t)f(t,u(t))=0, \,\, t \in (0,1)  \\
u(0)=0=u(1),
\end{cases}
  \end{equation}
with $q(t)=\left[ r_{2}\left(\frac{r_{1}}{r_{2}}\right)^{t} \log\left(\frac{r_{2}}{r_{1}}\right)\right]^{2}$ and $f(t,u)=h\left(r_{2}\left(\frac{r_{1}}{r_{2}}\right)^{t} ,u\right)$. On the other hand, if $N \geq 3$, then the transformation $t=-\frac{A}{r^{N-2}}$ and $u(t)=v(r)$ transforms the system (\ref{e1}) to the BVP (\ref{e4}), where $A=\frac{(r_{1}r_{2})^{N-2}}{r_{2}^{N-2}-r_{1}^{N-2}}$, $ q(t)=(N-2)^{-2}\frac{A^{2/(N-2)}}{(B-t)^{2(N-1)/(N-2)}}$, $B=\frac{r_{2}^{N-2}}{r_{2}^{N-2}-r_{1}^{N-2}}$  and $f(t,u)=h\left(\left(\frac{A}{B-t} \right)^{1/(N-2)}, u \right)$. The function $q(t)$ defined in (\ref{e4}) is well-defined, continuous and bounded between positive constants in the interval $[0,1]$.\\

\ \ \ \ Since we are interested in finding sufficient conditions for the existence of positive radial solutions of (\ref{e1}),
 it is equivalent to study the existence of positive solutions of (\ref{e4}).
 In this paper, we provide some new sufficient conditions for the existence of positive solutions of (\ref{e4}).\\

  \ \ \ \ In \cite{r18}, Iturriaga et. al used Krasnoselskii fixed point theorem for the existence of a positive solution of (\ref{e4}) for $\lambda$ small, and sub and super solution method for the existence of two positive solutions of (\ref{e4}) for $\lambda$ large. The main focus of the work in \cite{r18} is on the use of local superlinearity of the nonlinear function $f$ at $\infty$. Motivated by the work of Hai and Qian \cite{r15} for first order delay differential equations, and  the work of Gatica and Kim \cite{r12} for second order multipoint boundary value problems, we shall use two fixed point theorems by Gatika and Smith \cite{r13} to provide ranges on the parameter $\lambda$ in (\ref{e4})  to obtain sufficient conditions for the existence of positive solutions.\\

  \ \ \ \ Our main results are
  \begin{theorem}\label{thm1.1}
  Let
  \begin{equation}\label{e5}
f_{0}: \lim_{u \to 0+}\frac{f(t,u)}{u}=0 \,\,\,\, \text{uniformly \,\, in} \,\,  t \in (0,1).
  \end{equation}
  Then for each $R>0$, there exists a constant $\lambda_{R}>0$, $\lambda_{R}$ large enough such that for $\lambda>\lambda_{R}$, (\ref{e4}) has a positive solution $u(t)$ with $\sup_{t \in [0,1]} u(t) \leq R$.
  \end{theorem}

  \begin{theorem}\label{thm1.2}
  Let
  \begin{equation}\label{e6}
f_{0}: \lim_{u \to 0+}\frac{f(t,u)}{u}=\infty \,\,\,\, \text{uniformly \,\, in} \,\,  t \in (0,1).
  \end{equation}
  Then for each $R>0$, there exists a constant $\lambda_{R}>0$ such that for $\lambda<\lambda_{R}$, (\ref{e4}) has a positive solution $u(t)$ with $\sup_{t \in [0,1]} u(t) \leq R$.
  \end{theorem}

  \begin{theorem}\label{thm1.3}
  Let
  \begin{equation}\label{e7}
f_{\infty}: \lim_{u \to \infty}\frac{f(t,u)}{u}=0 \,\,\,\, \text{uniformly \,\, in} \,\,  t \in (0,1).
  \end{equation}
  Then for each $r>0$, there exists a constant $\lambda_{r}>0$ such that for $\lambda>\lambda_{r}$, (\ref{e4}) has a positive solution $u(t)$ with $\min_{t \in [0,1]} u(t) \geq r$.
  \end{theorem}

    \begin{theorem}\label{thm1.4}
Assume that there exist positive constants $0< \alpha < \beta <1$ such that
  \begin{equation}\label{e8}
f_{\infty}: \lim_{u \to \infty}\frac{f(t,u)}{u}=\infty \,\,\,\, \text{uniformly \,\, in} \,\,  t \in [\alpha,\beta].
  \end{equation}
  Then for each $r>0$, there exists a constant $\lambda_{r}>0$ such that for $\lambda<\lambda_{r}$, (\ref{e4}) has a positive solution $u(t)$ with $\min_{t \in [0,1]} u(t) \geq r$.
  \end{theorem}

\ \ \ \ Now we provide examples that strengthens our results.
\begin{example}\label{exa1.1}
Consider the elliptic equations in $\mathbb{R}^{3}$
\begin{equation}\label{ex1}
\begin{cases}
\Delta u+\lambda |x|^{4}\frac{u^{2}}{1+u}=0, \,\, t \in (0,1),\\
u(0)=0=u(1),
\end{cases}
\end{equation}
 Here $f(t,u)=\frac{u^{2}}{1+u}$. Clearly, $\lim_{u \to 0}\frac{f(t,u)}{u}=\lim_{u \to 0}\frac{u}{1+u}=0$ and $\lim_{u \to \infty}\frac{f(t,u)}{u}=\lim_{u \to \infty}\frac{u}{1+u}=1$ implies that Theorem \ref{thm1.1} can be applied to this example, where as Theorems \ref{thm1.2}--\ref{thm1.4} cannot be applied to this example. Now $m_{R}=\min_{u \in [0,R], t \in [1/4,3/4]}\frac{u}{1+u}<1$ and $\int^{3/4}_{1/4}G(s,s)\,ds=11/96$ implies that $\lambda_{R}>1536/11=139.636364$. By Theorem \ref{thm1.1}, the problem \ref{ex1} has a positive solution for $\lambda>1536/11$. Using Matlab, the value of $R$ is found to be $2.7679 \times 10^{-14}$. This is illustrated in  Figure 1a. 
\end{example}

\begin{example}\label{exa1.2}
Consider the elliptic equations in $\mathbb{R}^{3}$
\begin{equation}\label{ex2}
\begin{cases}
\Delta u+\lambda |x|^{4}(\sqrt{u}+\frac{u}{2})=0, \,\, t \in (0,1),\\
u(0)=0=u(1),
\end{cases}
\end{equation}
 Here $f(t,u)=\sqrt{u}+\frac{u}{2}$. Clearly, $\lim_{u \to 0}\frac{f(t,u)}{u}=\lim_{u \to 0}\frac{1}{u}(\sqrt{u}+\frac{u}{2})=\infty$ and $\lim_{u \to \infty}\frac{f(t,u)}{u}=\lim_{u \to \infty}\frac{1}{u}(\sqrt{u}+\frac{u}{2})=\frac{1}{2}$ implies that Theorem \ref{thm1.2} can be applied to this example, where as Theorems \ref{thm1.1}, \ref{thm1.3} and \ref{thm1.4} cannot be applied to this example. Now 
 \[M_{R}=\max_{u \in [0,R], t \in [1/4,3/4]}\frac{1}{u}(\sqrt{u}+\frac{u}{2})\geq\frac{1}{\sqrt{R}}+\frac{1}{2} \,\, \text{ and} \,\,  \int^{1}_{0}G(s,s)\,ds=1/6\] implies that $\lambda_{R}<\frac{12\sqrt{R}}{2+\sqrt{R}}\leq 12$. By Theorem \ref{thm1.2}, the problem \ref{ex2} has a positive solution for $\lambda<12$.  Using Matlab, the value of $R$ is found to be $5.1897$. This is illustrated in  Figure 1b.

\end{example}

\begin{example}\label{exa1.3}
Consider the elliptic equations in $\mathbb{R}^{3}$
\begin{equation}\label{ex3}
\begin{cases}
\Delta u+\lambda |x|^{4}\frac{u}{1+u}=0, \,\, t \in (0,1),\\
u(0)=0=u(1),
\end{cases}
\end{equation}
 Here $f(t,u)=\frac{u}{1+u}$. Clearly, $\lim_{u \to 0}\frac{f(t,u)}{u}=\lim_{u \to 0}\frac{1}{1+u}=1$ and $\lim_{u \to \infty}\frac{f(t,u)}{u}=\lim_{u \to \infty}\frac{1}{1+u}=1$ implies that Theorem \ref{thm1.3} can be applied to this example, where as Theorems \ref{thm1.1}, \ref{thm1.2} and \ref{thm1.4} cannot be applied to this example. Now $m_{r}=\min_{u \in [0,r], t \in [0,1]}\frac{1}{1+u}<1$ and $\int^{3/4}_{1/4}G(s,s)\,ds=11/96$ implies that $\lambda_{r}>384/11=34.9090909$. By Theorem \ref{thm1.3}, the problem \ref{ex1} has a positive solution for $\lambda>384/11$.  Using Matlab, the value of $r$ is found to be $0.0108$. This is illustrated in  Figure 1c.
\end{example}

\begin{example}\label{exa1.4}
Consider the elliptic equations in $\mathbb{R}^{3}$
\begin{equation}\label{ex4}
\begin{cases}
\Delta u+\lambda |x|^{4}(u^{3}+\frac{u}{2})=0, \,\, t \in (0,1),\\
u(0)=0=u(1), 
\end{cases}
\end{equation}
 Here $f(t,u)=(u^{3}+\frac{u}{2})$. Clearly, $\lim_{u \to \infty}\frac{f(t,u)}{u}=\lim_{u \to \infty}(u^{2}+\frac{1}{2})=\infty$ and $\lim_{u \to 0}\frac{f(t,u)}{u}=\lim_{u \to 0}(u^{2}+\frac{1}{2})=\frac{1}{2}$ implies that Theorem \ref{thm1.4} can be applied to this example, where as Theorems \ref{thm1.1}--\ref{thm1.3} cannot be applied to this example. Now $M_{R}=\max_{u \in [0,r], t \in [0,1]}(r^{2}+\frac{1}{2})$  implies that $\lambda_{r}<\frac{12}{2r^{2}+1}$. Setting $g(r)=\frac{12}{2r^{2}+1}$, we see that $g(r)$ attains its maximum 12 at $r=0$.  By Theorem \ref{thm1.4}, the problem \ref{ex2} has a positive solution for $\lambda<12$.  Using Matlab, the value of $r$ is found to be $8.2207 \times 10^{-10}$. This is illustrated in  Figure 1d.
\end{example}

\begin{figure}[h!]
	\centering
	\subfloat[Example 1.1]{\includegraphics[width=0.5\linewidth]{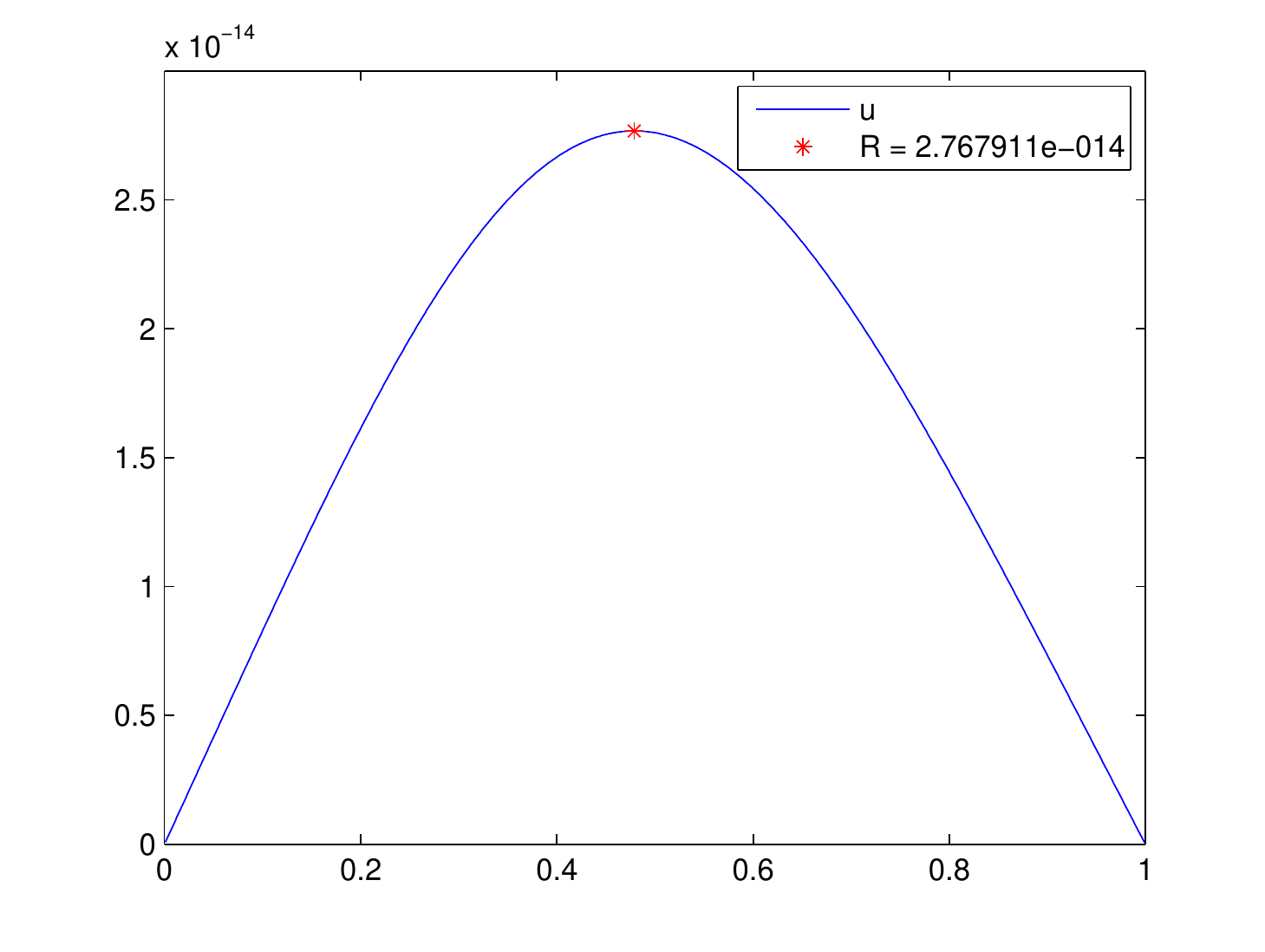}}
	\subfloat[Example 1.2]{\includegraphics[width=0.5\linewidth]{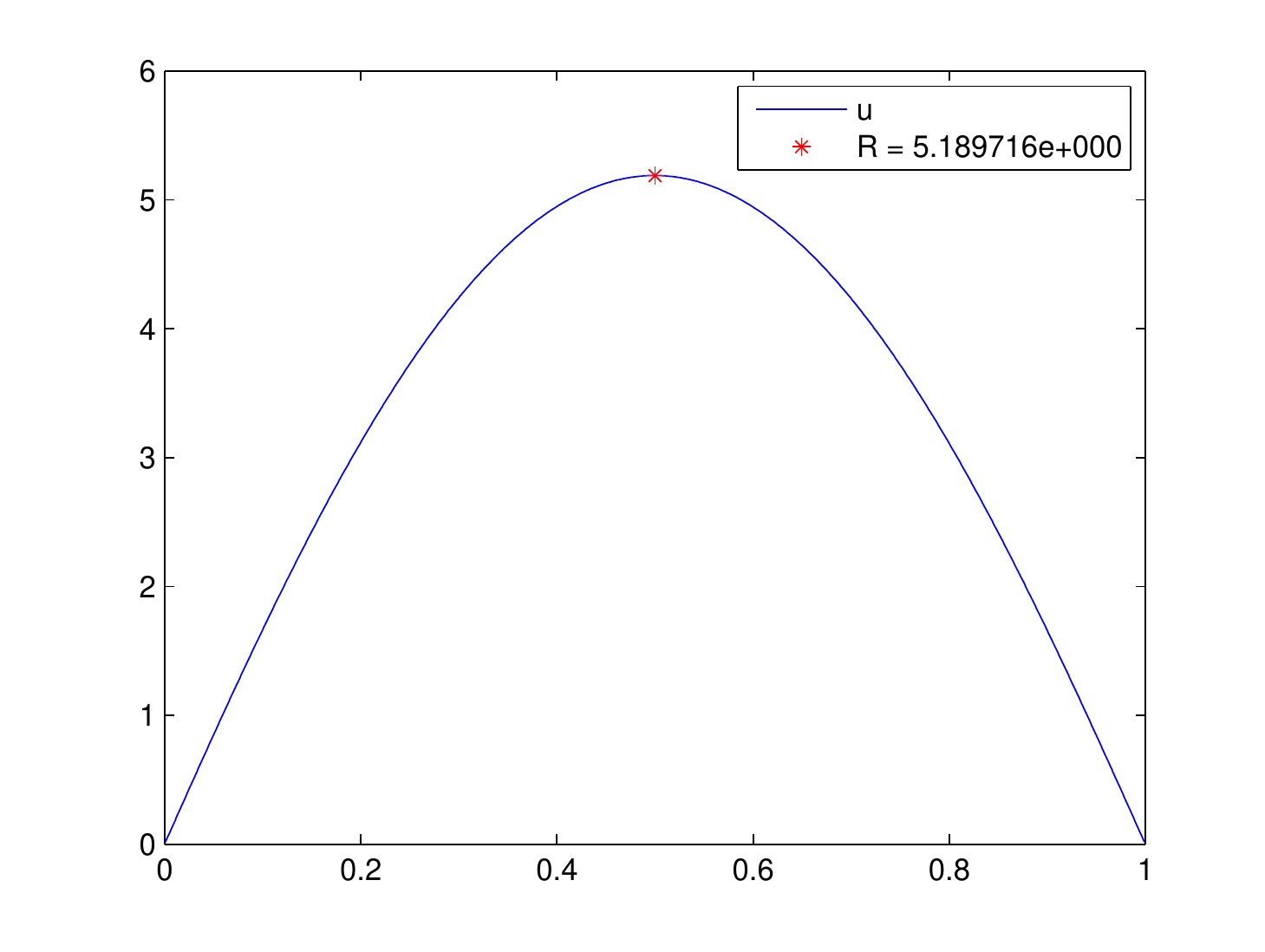}}\\
	\subfloat[Example 1.3]{\includegraphics[width=0.5\linewidth]{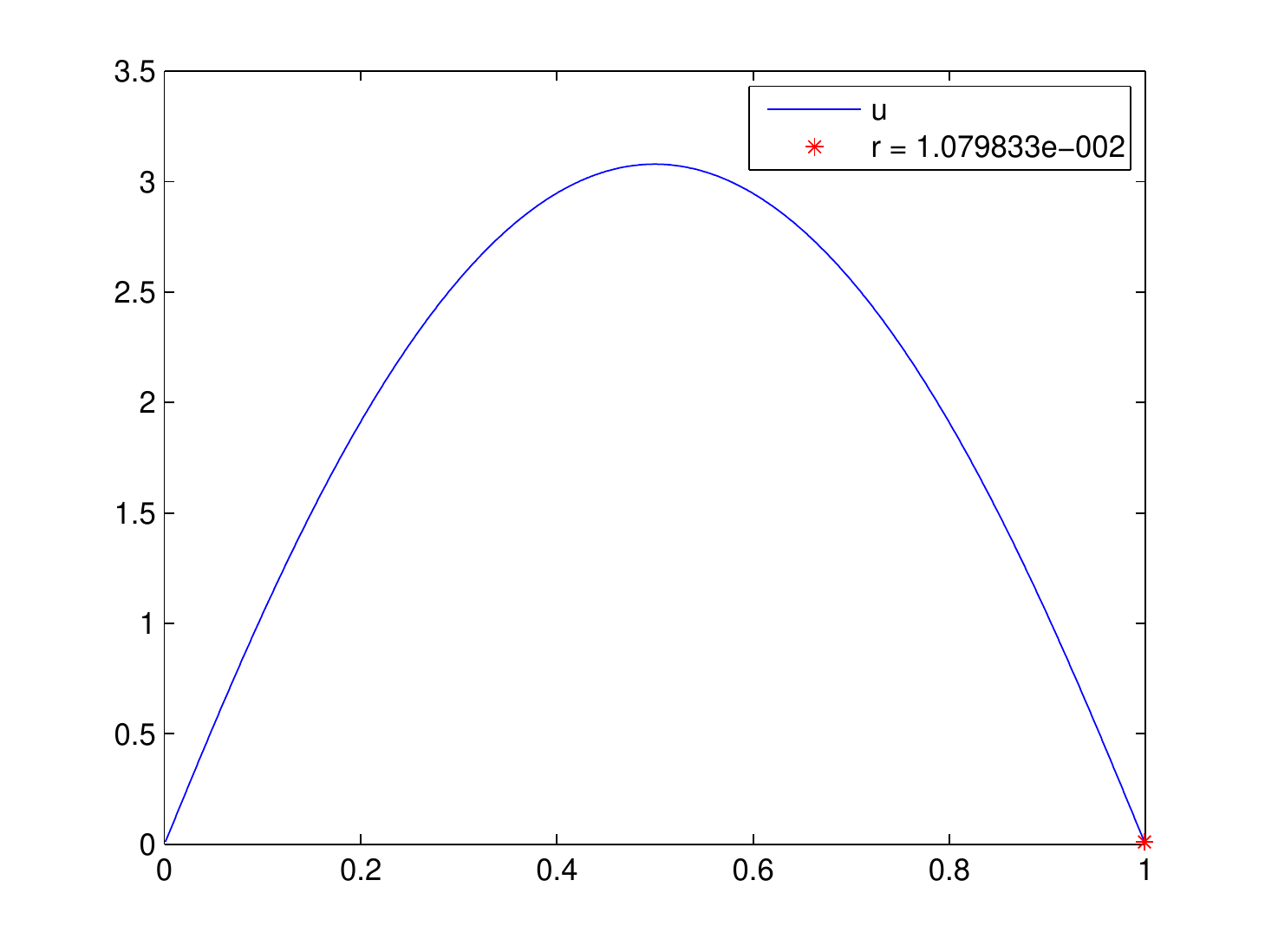}}
	\subfloat[Example 1.4]{\includegraphics[width=0.5\linewidth]{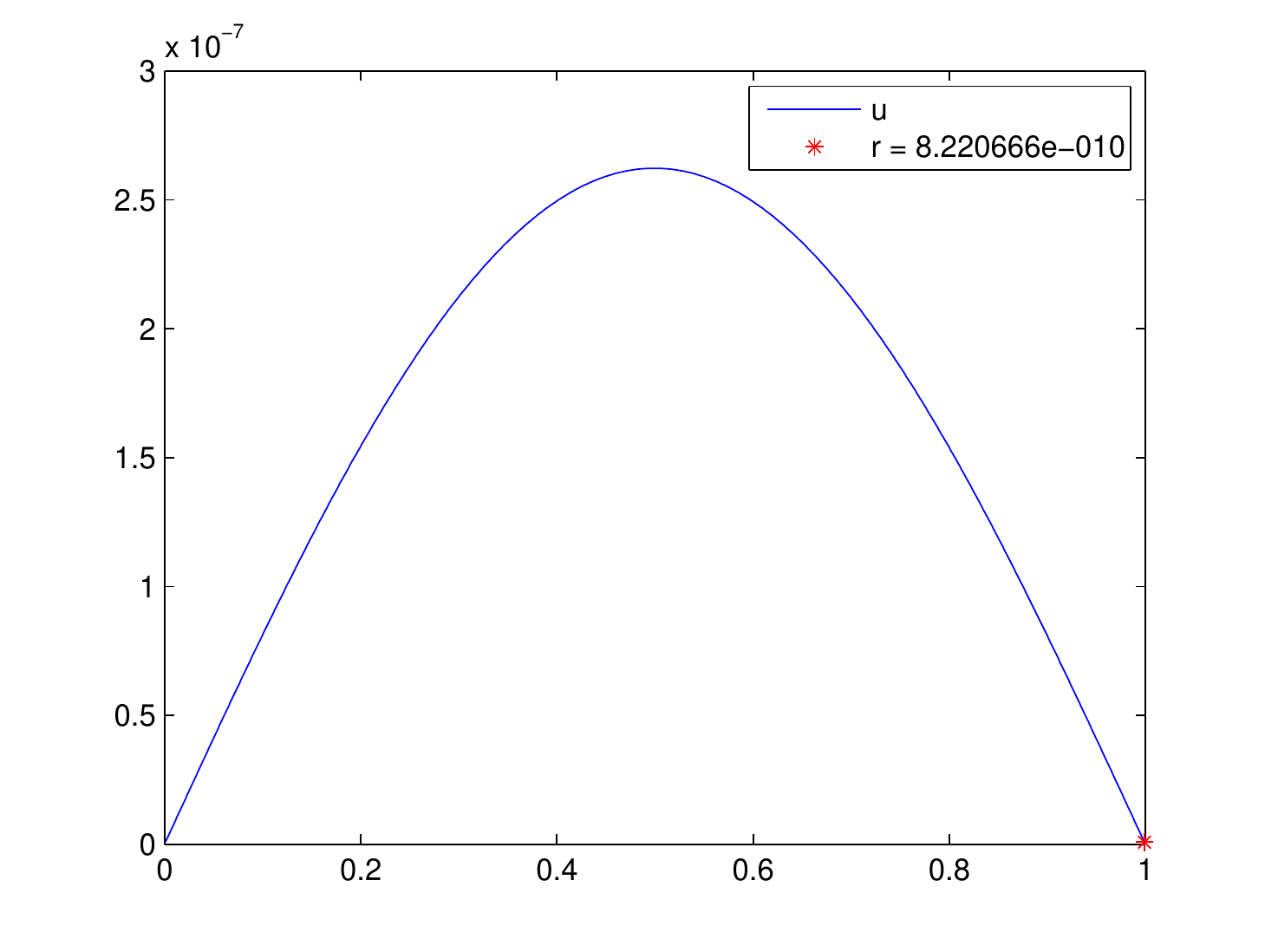}}
	\caption{}
\end{figure}

\ \ \ \   For our next two theorems, we  consider the eigen-value problem
\begin{equation}\label{e19}
\begin{cases}
-u''(t) = \lambda m(t)u(t), \, 0<t<1\\
u(0) = u(1) = 0,
\end{cases}
\end{equation}
where $m : [0,1] \to [0,\infty]$ is a continuous function. It is known that there exists a sequence of positive eigen values, which we denote by $\lambda_{i,m}(i>0),$ provided $m>0$ in a set of positive measure. In particular, $\lambda_{1,m}$ is called the first eigen value of the problem (\ref{e19}) and the associated eigen function, denoted by, $\phi_{1,m}$ satisfies the properties $\phi_{1,m} > 0, \, \phi_{1,m}(0) = \phi_{1,m}(1) = 0$ with $\phi'_{1,m}(0) >0$ and $\phi'_{1,m}(1)<0.$ The following two theorems provide inequalities on the function $f(t,u)$ and ranges on $\lambda$, using eigen values and their corresponding eigen functions for  the existence of positive solutions of (\ref{e4}).

\begin{theorem}\label{thm4.1}
Assume that there exist a continuous function $b:[0,1] \to [0,\infty]$ and positive constants $c$,$\delta$ and $R$ with $c>1$ and $0<\delta<R$ such that
\begin{description}
    \item[(H1)] $f(t,u)\leq b(t)u$ for $u \in (0,\delta)$ uniformly in $t \in (0,1)$
    \end{description}
  and
  \begin{description}
    \item[(H2)] $f(t,u) \geq cb(t)u$ for $u\geq R$ uniformly in $t \in (0,1).$
  \end{description}
holds. Then the BVP (\ref{e4}) has a positive solution for every
\[\frac{\lambda_{1 , qb}}{c} < \lambda < \lambda_{1qb}.\]
\end{theorem}

\begin{theorem}\label{thm4.2}
Assume that there exists a continuous function $b: [0,1] \to [0,\infty]$ and positive constants $c,\delta$ and $R$ with $c>1$ and $0<\delta<R$ such that
\begin{description}
    \item[(H3)] $f(t,u)\leq b(t)u$ for $u \geq R$ uniformly in $t \in (0,1)$
    \end{description}
and
 \begin{description}
    \item[(H4)] $f(t,u) \geq cb(t)u$ for $u\in (0,\delta)$ uniformly in $t \in (0,1).$
  \end{description}
  holds. Then the BVP (\ref{e4}) has a positive solution for every
  \[\frac{\lambda_{1,qb}}{c} < \lambda < \lambda_{1,qb}.\]
   \end{theorem}

\begin{example}\label{exa4.1}
 Consider the elliptic equations in $\mathbb{R}^{3}$
\begin{equation}\label{ex5}
\begin{cases}
\Delta u+\lambda |x|^{4}u^{3}=0, \,\, t \in (0,1),\\
u(0)=0=u(1), 
\end{cases}
\end{equation}
To the equation (\ref{ex5}), we associate the second order ODE
\begin{equation*}
\begin{cases}
u''(t)+\lambda u^{3}(t)=0, \, 0<t<1\\
u(0) = u(1) = 0.
\end{cases}
\end{equation*}
Here $f(t,u)=u^{3}$ and $q(t)=1$. Set $b(t)=1$. It is easy to verify that $\lambda_{1,qb} =\pi^{2}$ is the first eigen value of the equation 
\begin{equation*}
\begin{cases}
u''(t)+ \lambda u(t)=0, \, 0<t<1\\
u(0) = u(1) = 0.
\end{cases}
\end{equation*}
Let $c>1$ be a constant. Choose  $\delta \in (0,1)$ such that $c \delta^{2}>1$. Then for $u \in (0,\delta)$, we have $f(t,u) = u^{3}<u$. Hence the condition $(H1)$ of Theorem \ref{thm4.1} is satisfied.  Set $R=c\delta$. Then for $u \geq R$, we have $f(t,u)=u^{3}>u^{2}u>c^{2}\delta^{2}u=c.c\delta^{2}u>cu$. Thus, the condition $(H2)$ of Theorem \ref{thm4.1} is satisfied. Hence, by Theorem \ref{thm4.1}, (\ref{ex5}) has a positive solution for every $\frac{\pi^{2}}{c}<\lambda <\pi^{2}$. If we set $\lambda=\pi^{2}/2$, then the maximun amd minimum value of the solution $u(t)$, using MATLAB, is  $2.1278 \times 10^{-36}$ and $4.6998 \times 10^{-39}$.
\end{example}

\begin{example}
Consider two constants  $\delta \in (0,1)$ and $R>1$.  Choose $c=\frac{1}{\delta^{2}}$. Then by Theorem \ref{thm4.2}, the equation
\begin{equation*}
\begin{cases}
\Delta u+\lambda |x|^{4}u^{-1}=0, \,\, t \in (0,1),\\
u(0)=0=u(1), 
\end{cases}
\end{equation*} 
has a positive solution $u(t)$ for each $\lambda \in (\delta^{2}\pi^{2},\pi^{2})$.
\end{example}

\ \ \ \ In \cite{r17}, Henderson and Wang provide ranges on $\lambda$, depending on $f_{0}$ and $f_{\infty}$, to obtain at least one positive solution of (\ref{e4}), provided that both $f_{0}$ and $f_{\infty}$ exist (see Theorem 2 and Theorem 3 in \cite{r17}). In a similar way, Lan and Webb \cite{r20} proved several existence theorem on the positive solution of (\ref{e4}) within some particular range on $\lambda$, depending upon $f_{0}$ and $f_{\infty}$. The basic idea of the proofs in \cite{r17} and \cite{r20} are Krasnoselskii's fixed point theorem and a fixed point theorem due to Amann. \\

\ \ \ \ Results similar to Theorems \ref{thm1.1}--\ref{thm1.4} can be found in \cite{r32}. Wang used fixed point index approach to obtain positive solutions of a system of equations, see Theorem 1.2 (a) and (b) in \cite{r32}. Wang \cite{r33} obtained an existence of a positive solutions of (\ref{e4}) under the assumption that $f$ is sublinear.  Theorems \ref{thm4.1}--\ref{thm4.2} are based on inequalities to have positive solutions of the BVP (\ref{e4}), which are completely new in the literature. The ranges on $\lambda$ in Theorems \ref{thm4.1}--\ref{thm4.2} are completely dependent on the first eigen value of the eigen value problem (\ref{e19}) with $m(t)=q(t)b(t)$.

\ \ \ \ This work has been divided into three sections. Section 1 is Introduction. We provide the statements of our theorems.  In Section 2, we provide some basic results of this paper. The proof of the  Theorems \ref{thm1.1}--\ref{thm4.2} are given in Section 3.

\section{Preliminaries}

\ \ \ \ We consider the Banach space $X = C([0,1])$  endowed with the norm
\begin{equation}\label{e9''}
\| x\|=\max_{t \in [0,1]}|x(t)|,
\end{equation}
and a cone $K$ on $X$ by
\begin{equation}\label{e9'}
    K=\{ u \in X; u(t) \geq 0, t \in (0,1), u(0)=0=u(1)    \}.
\end{equation}
Define an operator $T:K \to X$ by
\begin{equation}\label{e9}
    (Tu)(t)=\lambda \int^{1}_{0}G(t,s)q(s)f(s,u(s))\,ds,
\end{equation}
where $G(t,s)$ is the Green's function in the interval $(0,1)$, given by
 \begin{equation}\label{e10}
G(t,s)=\begin{cases}
 s(1-t); \,\, 0 \leq s \leq t \leq 1  \\
 t(1-s); \,\, 0 \leq t \leq s \leq 1.
\end{cases}
  \end{equation}

It is proved in \cite{r9} that
\begin{align}
G(t,s) > 0  \,\, & \text{on} \,\, (0,1) \times (0,1) \nonumber \\
 G(t,s) \leq G(s,s) = &  s(1-s), \,\, 0 \leq s,  t \leq 1, \nonumber
  \end{align}
  and\\
 \[  G(t,s) \geq \frac{1}{4}G(s,s) = \frac{1}{4}s(1-s), \,\, \frac{1}{4} \leq t \leq \frac{3}{4}, \,\, 0 \leq s \leq 1. \]

\ \ \ \ We shall use the following fixed point results in a cone \cite{r12}, which are the revised version of theorems due to Gustafson and Schmitt \cite{r13'}.
 \begin{theorem}\label{thm2.1}
  Let $X$ be a Banach space and $K$ be a cone in $X$. Let $r$ and $R$ be real numbers with $0<r<R$,
  \ \[ D = \{u \in K ; r \leq \|u\| \leq R\}, \]
  and Let $T: D \to K$ be a compact continuous operator such that
  \begin{description}
    \item[(a)] $u \in D, \mu < 1, u = \mu T u\Longrightarrow \|u\| \neq R$
    \item[(b)] $u \in D, \mu > 1, u = \mu T u\Longrightarrow \|u\| \neq r$
    \item[(c)] $\inf\limits_{\|u\|=r}\|Tu\| \neq 0.$
  \end{description}
   Then $T$ has a fixed point in $D$.
   \end{theorem}

   \begin{theorem}\label{thm2.2}
  Let $X$ be a Banach space and $K$ be a cone in $X$. Let $r$ and $R$ be real numbers with $0<r<R$,
  \ \[ D = \{u \in K ; r \leq \|u\| \leq R\}, \]
  and Let $T: D \to K$ be a compact continuous operator such that
  \begin{description}
    \item[(a)] $u \in D, \mu > 1, u = \mu T u\Longrightarrow \|u\| \neq R$
    \item[(b)] $u \in D, \mu < 1, u = \mu T u\Longrightarrow \|u\| \neq r$
    \item[(c)] $\inf\limits_{\|u\|=R}\|Tu\| > 0.$
  \end{description}
  Then $T$ has a fixed point in $D$.
   \end{theorem}

\ \ \ \ In order to satisfy the condition (c) in Theorem 2.1 and Theorem 2.2, we shall make an extensive use of the following lemma, given in  \cite{r12}.
\begin{lemma}\label{Lemma 2.3}
 Let $\phi : [0,1] \to [0,\infty)$ be a continuous function whose graph is concave down, and let $ \|\phi\| : max\{ \phi(t) : t \in [0,1]\}.$ Then, for any $t \in [\alpha,1-\alpha]$ with $0<\alpha<\frac{1}{2}$, we have $\alpha\|\phi\| \leq \phi(t)$.
 \end{lemma}

 \begin{lemma}\label{Lemma 2.4}
 If $R > 0$ is a real number, then\\
\[ \inf  \{\|Tu\| ; u\in K \,\,  and  \,\,   \|u\| = R\} > 0 \]
 for any solution $u$ of (\ref{e4}).
  \end{lemma}
\noindent \textbf{\underline{Proof}:}
Clearly, $u(t)$ in a solution of (\ref{e4}) if and only if $Tu = u$. Since $(Tu)'' = -\lambda q(t)f(t,u)$, then the graph of $Tu$ is always concave down, and the graph of $u$ is concave down. Hence, for $\theta \in (0,\frac{1}{2})$, it follows from lemma 2.1 that
\begin{align}
  u(t) \geq \theta \|u\| \,\,& \text{for}\,\, t \in (\theta,1-\theta) \nonumber
\end{align}
Let $\theta = \frac{1}{4}.$ Since $f(t,u) > 0$ for $t \in [\frac{1}{4},\frac{3}{4}]$ and $ u \in [\frac{1}{4}R,R],$ then for
\[p = \inf \{f(t,u);(t,u) \in [\frac{1}{4},\frac{3}{4}] \times [\frac{R}{4},R]\} > 0\]
and\\
\[q = \inf\limits_{t \in [\frac{1}{4},\frac{3}{4}]} t(1-t)q(t),\]
we have
\begin{align}
(Tu)(t) = & \lambda\int^{1}_{0}G(t,s)q(s)f(s,u(s))\,ds, \nonumber \\
& > \lambda\int^{1}_{0}G(t,s)q(s)f(s,u(s))\,ds \nonumber \\
& \geq \lambda \frac{1}{4}\int^{\frac{3}{4}}_{\frac{1}{4}}s(1-s)q(s)f(s,u(s))\,ds \nonumber \\
& \geq \lambda \frac{1}{4}  pq(\frac{3}{4}-\frac{1}{4}) = \frac{\lambda}{8}pq > 0, \nonumber \
\end{align}
and so $\|Tu\| \geq \frac{\lambda}{8}pq > 0$ for all $u \in K$ with $\|u\| = R.$ The lemma is proved.

\section{\underline{Proof of the Main Results:}}
\ \ \ \ In this section, we consider the operator $T$ defined in (\ref{e9}), and the Banach space X in (\ref{e9''}) and cone $K$ in (\ref{e9'}).\\

\noindent \textbf{\underline{Proof of Theorem \ref{thm1.1}}:}
Let $R > 0$. Choose $\lambda_{R}>0, \, \lambda_{R}$ large enough such that
\[\lambda_{R} > \frac{16}{m_{R}\int^{\frac{3}{4}}_{\frac{1}{4}}G(s,s)q(s)ds},\]
where\\
\begin{equation}\label{e11}
 m_{R} = \min_{t \in [\frac{1}{4},\frac{3}{4}]_{0 \leq u \leq R}} \frac{f(t,u)}{u}.
\end{equation}\\
Let $\lambda > \lambda_{R}.$ By (\ref{e5}), there exists $r \in (0,R)$ and $\epsilon > 0$ such that $f(t,u) \leq \epsilon u$ for $0 < u \leq r$ and
\[0 < \epsilon < \frac{1}{\lambda \int^{1}_{0}G(s,s)q(s)ds}.\]
Consider
\[ D = \{u \in K ; r\leq u(t)\leq R, t \in (0,1)\} \]
Using Arzela - Ascoli lemma, we can prove that $T: D \to K$ is compact and continuous. In order to complete the proof of the  theorem, we shall use Theorem 2.2.

\ \ \ \ Let $u \in D$ be such that $u = \mu Tu$ and $\mu > 1$, that is
\begin{equation}\label{e12}
u(t) = \mu \lambda \int^{1}_{0} G(t,s)q(s)f(s,u(s))ds, \,\, \mu>1
\end{equation}
We claim that (\ref{e12}) has no solution with $\|u\| = R.$ Suppose that (\ref{e12}) has a solution $u_{0} (t)$ with $\|u_{0}\| = R.$ Without any loss of generality, we assume that $u_{0} (t)\geq 0$ for $t \in (0,1).$ Then
\begin{align}
\|u_{0}\|\geq\min\limits_{t \in [\frac{1}{4},\frac{3}{4}]} u_{0} (t) & \geq \mu \lambda \int^{1}_{0}(\min\limits_{t \in [\frac{1}{4},\frac{3}{4}]} G(t,s)) q(s) f(s,u_{0}(s))ds \nonumber \\
&\geq \frac{\mu \lambda}{4} \int^{1}_{0} G(s,s) q(s) f(s,u_{0}(s))ds \nonumber\\
&\geq \frac{\mu \lambda}{4} \int^{\frac{3}{4}}_{\frac{1}{4}} G(s,s) q(s) f(s,u_{0}(s))ds \nonumber \\
&\geq \frac{\mu \lambda m_{R}}{4} \int^{\frac{3}{4}}_{\frac{1}{4}} G(s,s) q(s) u_{0}(s)ds. \label{e13}
\end{align}
Since $u_{0}(t)$ is a solution of (\ref{e12}), then it satisfies
\begin{equation}\label{e14}
\begin{cases}
 u_{0}''(t) = - \lambda \mu q(t)f(t,u_{0}(t)),  \\
 u_{0}(0) = u_{0} (1) = 0.
\end{cases}
  \end{equation}
Since the graph of $u_{0} (t)$ is concave down, then by Lemma 2.3 with $\alpha = \frac{1}{4}$, we have
\begin{equation}\label{e15}
u_{0}(t) \geq \frac{1}{4} \|u_{0}\|
\end{equation}
Using (\ref{e15}) in (\ref{e13}), we have\\
\begin{align}
R = \|u_{0}\| & \geq \frac{\lambda \mu m_{R}}{16} \|u_{0}\| \int^{\frac{3}{4}}_{\frac{1}{4}} G(s,s)q(s)ds \nonumber \\
&\geq \lambda_{R}\mu \frac{m_{R}}{16} \|u_{0}\| \int^{\frac{3}{4}}_{\frac{1}{4}} G(s,s)q(s)ds \nonumber \\
& > \mu  \|u_{0}\| = \mu R > R, \nonumber
\end{align}
a contradiction. Hence our claim holds, that is, $\|u\| = R.$

\ \ \ \  Next, let $u \in D$ with $u=\mu Tu$ for some $\mu \in (0,1).$ We claim that $\|u\| \neq r$. Suppose that $\|u\| = r$.  Then
 \begin{align}
 r = \|u\| & \leq \mu \lambda \int^{1}_{0}G(t,s)q(s)f(s,u(s))ds \nonumber\\
 & \leq \mu \lambda \epsilon \|u\| \int^{1}_{0}G(s,s)q(s)ds \nonumber \\
  & \leq \mu \lambda \epsilon r \int^{1}_{0}G(s,s)q(s)ds \nonumber\\
  & < \mu r < r, \nonumber
 \end{align}
 a contradiction. Hence $\|u\| \neq r.$

\ \ \ \  By Lemma \ref{Lemma 2.4} and Theorem \ref{thm2.2}, the BVP (\ref{e4}) has a positive solution $u$ in $D$ with  $\sup\limits_{t \in [0,1]} u(t) \leq R$. This completes the proof of the theorem.\\

\noindent \textbf{\underline{Proof of Theorem \ref{thm1.2}}:}
By (\ref{e6}), there exists a constant $r \in (0,R)$ and $B$ with
\[ B > \frac{16}{\lambda}(\int^{\frac{3}{4}}_{\frac{1}{4}}G(s,s)q(s)ds)^{-1}\]
such that
\[ f(t,u) > Bu \,\, \text{ for} \,\, u \in (0,r], \,  0<t<1\]
Now, we consider the set
\[D = \{u \in K : \,  r \leq u \leq R, \, t \in (0,1)\} \]
where $K$ is the cone given in (\ref{e9'}). We consider  the operator $T$ on $X$ as in (\ref{e9}). An application of Arzela - Ascoli lemma proves that $T : D \to K$ is compact and continuous. Choose $\lambda _{R} > 0$ small enough such that
\[\lambda_{R} \leq \frac{1}{M _{R}\int^{1}_{0}G(s,s)q(s)ds},\]
where
\begin{equation}\label{e16}
M_{R}= \max\limits_{0 \leq u \leq R, \, 0 \leq t \leq 1} \frac{f(t,u)}{u}.
\end{equation}
We shall use Theorem 2.1 to prove the  Theorem. Let $u \in D$ be such that $u = \mu Tu$ for some $\mu \in (0,1)$. In this case, we claim that $\|u\| \neq R$. On the contrary, suppose that $\|u\| = R$. Then
\begin{align}
u(t) =  & \mu \lambda \int^{1}_{0} G(t,s)q(s)f(s,u(s))ds \nonumber \\
& < \mu \lambda_{R} \int^{1}_{0} G(s,s)q(s)f(s,u(s))ds \nonumber \\
& <\mu \lambda_{R} M_{R} \int^{1}_{0} G(s,s)q(s)u(s)ds \nonumber
\end{align}
implies that
\[R = \|u(t)\| < \mu \lambda_{R}.\|u\|.M_{R}\int^{1}_{0}G(s,s)q(s)ds < \mu R < R,\]
a contradiction. Hence $\|u\| \neq R.$

 \ \ \ \ Next, suppose that $u \in D$ and $u = \mu Tu$ for some $\mu > 1.$ We claim that $\|u\| \neq r.$ If possible, suppose that $\|u\| = r.$ Since $u = \mu Tu$, then we have
 \begin{align}
 u(t) \geq  & \mu \lambda B\int^{1}_{0}G(t,s)q(s)u(s)ds \nonumber \\
  \geq & \mu \lambda B\int^{\frac{3}{4}}_{\frac{1}{4}}G(t,s)q(s)u(s)ds. \nonumber
 \end{align}
Hence
\begin{align}
\|u\| \geq \min\limits_{t \in [\frac{1}{4},\frac{3}{4}]} u(t) & \geq \mu \lambda B\int^{\frac{3}{4}}_{\frac{1}{4}}(\min\limits_{t \in [\frac{1}{4},\frac{3}{4}]}G(t,s))q(s)u(s)ds \nonumber \\
& \geq \frac{\lambda \mu B}{4} \int^{\frac{3}{4}}_{\frac{1}{4}}G(s,s)q(s)u(s)ds. \label{e17}
\end{align}
Since $ u =\mu Tu,$ then $u$ satisfies (\ref{e14}) and hence $u(t)$ satisfies the property (\ref{e15}). Consequently, we obtain
\begin{align}
 r = \|u\| & \geq \frac{\mu \lambda B}{4}\frac{1}{4}\|u\| \int^{\frac{3}{4}}_{\frac{1}{4}}G(s,s)q(s)ds \nonumber\\
&\geq \mu . B. \frac{\lambda}{16} r \int^{\frac{3}{4}}_{\frac{1}{4}}G(s,s)q(s)ds \nonumber \\
& > \mu r > r, \nonumber
\end{align}
which is a contradiction. Hence our claim holds, that is, $\|u\| \neq r.$ Hence by Lemma \ref{Lemma 2.4} and Theorem \ref{thm2.1}, the BVP (\ref{e4}) has a positive solution $u(t)$ satisfying $r \leq u(t) \leq R,$ $t \in (0,1)$.  The theorem is proved.\\

\noindent \textbf{\underline{Proof of Theorem \ref{thm1.3}}:}
Let $r > 0$ and choose $\lambda_{r} > 0$ such that
\[\lambda _r > \frac{4}{m_r . \int^{\frac{3}{4}}_{\frac{1}{4}}G(s,s)q(s)ds},\]
where
\[m_r = \min\limits_{t \in [\frac{1}{4},\frac{3}{4}],  0 \leq u \leq r}\frac{f(t,u)}{u}.\]
By (\ref{e7}), we can find a constant $\epsilon > 0$ with
\[ \epsilon < \frac{1}{\lambda}\left(\int^{1}_{0}G(s,s)q(s)ds\right)^{-1}\]
and a constant $R_{0} > r$ such that $f(t,u) < \epsilon u $ for $u \geq R_{0}$.

\ \ \ \ We shall use Theorem \ref{thm2.1} to prove the  theorem. We claim that the equation $u = \mu Tu,\, 0<\mu <1  $
has no solution of norm $R,\, R\geq R_{0}.$ On the contrary, assume that there exists a sequence ${\{R_{n}\}}^{\infty}_{n=1}, \, R_{n} \to \infty$ as $ n \to \infty,\,R_{n} \geq R_{0}, \, n=1,2..,$ and a sequence ${\{\mu _n\}}^{\infty}_{n=1}$ of real numbers with $0<\mu_{n} < 1,$ and a sequence of functions ${\{u_{n}\}}^{\infty}_{n=1}$ with $\|u_{n}\| = R_{n}$ and
\begin{equation}\label{e18}
u_{n} = \mu _{n}Tu_{n} \, ,\, n=1,2,3 \cdots .
\end{equation}
Let $\{t_{n}\}$ be the unique point in $[0,1]$ such that $u_{n}(t_{n}) = \|u_{n}\|.$ Then from (\ref{e18}), we have
\begin{align}
R_{n} = u(t_{n}) = & \mu _{n} \lambda \int^{1}_{0} G(t_{n},s)q(s)f(s,u_{n}(s))ds \nonumber\\
 & \leq \mu _{n} \lambda \int^{1}_{0} G(s,s)q(s)f(s,u_{n}(s))ds \nonumber \\
 & \leq  \mu _{n} \lambda \epsilon \|u_{n}\|. \int^{1}_{0} G(s,s)q(s)ds \nonumber \\
 & < \mu _{\lambda} R_{n} < R_{n}, \nonumber
\end{align}
a contradiction. Hence our claim holds. Let us fix a real number $R > R_{0} \, .$ Then, by the above argument, we have that $u = \mu T u \, , \, 0<\mu <1$ has no solution with $\|u\| = R.$   Thus, if  we consider the set
\[D = \{u \in K; \, r\leq \|u\| \leq R, \, \, 0\leq t \leq 1\},\]
then, for the above choice of $R$, the condition (a) of Theorem \ref{thm2.1} is satisfied.

\ \ \ \ Now, we prove the condition (b) of Theorem 2.1 .  Let $u \in D$ with $u = \mu Tu, \, \mu >1.$ We claim that $\|u\| \neq r,$ If possible, let $u_{0} \in D$ be a solution of $u = \mu Tu, \, \mu >1.$ such that $\|u_{0}\| = r.$ Then
\[u_{0} (t) = \mu \lambda \int^{1}_{0} G(t,s)q(s)f(s,u_{0}(s))ds.\]
Since $u_{0} \in D$ with $u_{0} \in K$ and $0\leq u_{0}(t) \leq r$ with $\|u_{0}\| = r,$ then
\begin{align}
\min\limits_{t \in [\frac{1}{4},\frac{3}{4}]} u_{0}(t) = & \mu \lambda \int^{1}_{0}(\min \limits_{t \in [\frac{1}{4},\frac{3}{4}]}G(t,s))q(s)f(s,u_{0}(s))ds \nonumber \\
& \geq \mu \lambda_{r}\int^{1}_{0}(\min \limits_{t \in [\frac{1}{4},\frac{3}{4}]}G(t,s))q(s)f(s,u_{0}(s))ds \nonumber \\
& > \mu \lambda_{r}\int^{\frac{3}{4}}_{\frac{1}{4}}(\min \limits_{t \in [\frac{1}{4},\frac{3}{4}]}G(t,s))q(s)f(s,u_{0}(s))ds \nonumber \\
& >\frac{\mu}{4} \lambda_{r}\int^{\frac{3}{4}}_{\frac{1}{4}} G(s,s)q(s)f(s,u_{0}(s))ds \nonumber \\
& >\frac{\mu}{4} \lambda_{r} m_{r}\int^{\frac{3}{4}}_{\frac{1}{4}} G(s,s)q(s)u_{0}(s)ds \nonumber
\end{align}
Let $t_{0} \in [\frac{1}{4},\frac{3}{4}]$ be such that
\[u_{0}(t_{0}) = \min \limits_{t \in [\frac{1}{4},\frac{3}{4}]} u_{0}(t).\]
Then
\begin{align}
u_{0}(t_{0}) > & \frac{\mu}{4} \lambda_{r} m_{r} u_{0}(t_{0})\int^{\frac{3}{4}}_{\frac{1}{4}} G(s,s)q(s)ds \nonumber \\
&> \mu u_{0}(t_{0}) > u_{0}(t_{0}), \nonumber
\end{align}
a contradiction. Hence the condition (a) of Theorem \ref{thm2.1} is satisfied. The condition (c) of Theorem \ref{thm2.1} follows from Lemma \ref{Lemma 2.4} . By Theorem \ref{thm2.1}, BVP (\ref{e4}) has a positive solution $u(t)$ in $D$ satisfying $0 < r < \|u\| < R.$ This completes the proof of the theorem.\\

\textbf{\underline{Proof of Theorem \ref{thm1.4}}:}
Let $r>0$ be a constant. Choose $\lambda _r >0$ such that
\[\lambda_{r} < \frac{1}{M_{r} \int^{1}_{0}G(s,s)q(s)ds},\]
where
\[M_{r} = \max\limits_{0 \leq u \leq r,  0\leq t \leq 1} \frac{f(t,u)}{u} , \, \,\]
By (\ref{e8}), there exist constants $R_{0} > 0$ and $B > 0$ such that
\[f(t,u) \geq Bu  \, \, for \, \ u \geq R_{0},\]
where $B$ satisfies
\[\lambda B \int^{\frac{3}{4}}_{\frac{1}{4}}G(s,s)q(s)ds > 4.\]
We shall use Theorem \ref{thm2.2} to prove our theorem. We claim that for any $R > R_{0},$ the problem $u = \mu Tu, \, \mu >1$ has no solution with $\|u\| = R.$ If this is not true, then there exists a sequence $\{R_{n}\}^{\infty}_{n=1}, \, R_{n} \to \infty$\, as \,$ n \to \infty, \, R_{n} > R_{0}$ and a sequence $\{\mu_{n}\}^{\infty}_{n=1}$ of real numbers with $\mu_{n} > 1$ and a sequence of functions $\{u_{n}\}^{\infty}_{n = 1} $ with $\|u_{n}\| = R_{n}$ such that $u_{n}$ satisfies (\ref{e18}). Then we have
\begin{align}
u_{n}(t) = & \mu_{n}\lambda \int^{1}_{0}G(t,s)q(s)f(s,u_{n}(s))ds \nonumber \\
& \geq \mu_{n}B\lambda \int^{1}_{0}G(t,s)q(s)u_{n}(s)ds. \nonumber
\end{align}
Let $t^{*} \in [\frac{1}{4},\frac{3}{4}]$ be such that
\[\min \limits_{t \in [\frac{1}{4},\frac{3}{4}]} u_{n}(t) = u_{n}(t^{*}), \,\, n=1,2,\cdots.\]
Then
\begin{align}
u_{n}(t^{*}) = \min \limits_{t \in [\frac{1}{4},\frac{3}{4}]}u_{n}(t) & \geq \mu_{n}B\lambda \int^{1}_{0}(\min \limits_{t \in [\frac{1}{4},\frac{3}{4}]}G(t,s))q(s)u_{n}(s)ds \nonumber \\
& >  \mu_n B\lambda \int^{\frac{3}{4}}_{\frac{1}{4}}(\min \limits_{t \in [\frac{1}{4},\frac{3}{4}]}G(t,s))q(s)u_{n}(s)ds \nonumber \\
& \geq \frac{\mu_{n}B\lambda}{4}\int^{\frac{3}{4}}_{\frac{1}{4}}G(s,s)q(s)u_{n}(s)ds \nonumber \\
& \geq \frac{\mu_{n}B\lambda}{4}u_{n}(t^{*})\int^{\frac{3}{4}}_{\frac{1}{4}}G(s,s)q(s)ds \nonumber\\
& > \mu_{n} u_{n}(t^{*}) > u_{n}(t^{*}), \nonumber
\end{align}
a contradiction. hence, our claim holds. Fix $R \geq R_{0}$. Then for any $u \in K$ with $u = \mu Tu$ and $\mu>1$, we have $\|u\| \neq R$. Thus, if we consider the set
\[D \{ u \in K: \, r \leq \|u\| \leq R\},\]
then for the above choice of $R$, the condition (a) of Theorem \ref{thm2.2} is satisfied.

\ \ \ \  Let  $u \in D$ be such that $u=\mu Tu$ and $0 < \mu <1$. We claim that $\|u\| \neq r$. If possible, suppose that $\|u\|=r$. Then
\begin{align}
r = \|u\| & \leq \mu \lambda \int^{1}_{0}G(s,s)q(s)f(s,u(s))\,ds \nonumber \\
& \leq \mu \lambda M_{r} \int^{1}_{0}G(s,s)q(s)u(s)\,ds \nonumber \\
& \leq \mu \lambda_{r} M_{r} \|u\| \int^{1}_{0}G(s,s)q(s)\,ds \nonumber \\
& < \mu \|u\| < \|u\|=r, \nonumber
\end{align}
a contradiction. Hence, the condition (b) of Theorem \ref{thm2.2} is satisfied. The proof of condition (c) of Theorem \ref{thm2.2} is similar to the proof of Lemma \ref{Lemma 2.4}. By Theorem \ref{thm2.2}, the BVP (\ref{e4}) has a positive solution $u(t)$ with $\min_{t \in [0,1]}u(t) \geq r$. This completes the proof of the theorem.\\

\noindent \textbf{\underline{Proof of Theorem \ref{thm4.1}}:}
We shall use Theorem \ref{thm2.2}  to prove the theorem, Let $r\in(0,\delta).$ We claim that the integral equation
\begin{equation}\label{e20}
u(t) = \mu Tu, \, \, 0<\mu<1
\end{equation}
has no solution with norm $r.$ If possible, suppose that $u_{0}(t)$ is a solution of (\ref{e2}) with $\|u_{0}\| = r.$ Then $u_{0}(t)$ is a solution of the boundary value problem.
\begin{equation}\label{e21}
\begin{cases}
u''_{0}(t) + \lambda\mu q(t)f(t,u_{0}(t)) = 0, \, 0<\mu<1, \, 0<t<1\\
u_{0}(0) = 0 = u_{0}(1).
\end{cases}
\end{equation}
Multiplying (\ref{e21}) by $\phi_{1,qb}(t)$ and integrating both side from $0$ to $1$, we obtain
\begin{align}
-\int^{1}_{0}u''_{0}(t)\phi_{1,qb}(t)dt = &  \mu \lambda \int^{1}_{0}q(t)\phi_{1,qb}(t)f(t,u_{0}(t))dt \nonumber \\
& \leq \mu \lambda \int^{1}_{0}q(t)\phi_{1,qb}(t)b(t)u_{0}(t)dt. \label{e22}
\end{align}
Now,
\begin{align}
-\int^{1}_{0}u''_{0}(t)\phi_{1,qb}(t)dt = & \int^{1}_{0}u'_{0}(t)\phi'_{1,qb}(t)dt  \nonumber \\
= &  -\int^{1}_{0}u_{0}(t)\phi''_{1,qb}(t)dt \nonumber \\
= & \lambda_{1,qb}\int^{1}_{0}q(t)b(t)u_{0}(t)\phi_{1,qb}(t)dt \nonumber
\end{align}
implies, using  (\ref{e22}),  that
\begin{align}
\lambda_{1,qb}\int^{1}_{0}q(t)b(t)\phi_{1,qb}(t)u_{0}(t)dt & \leq \mu \lambda \int^{1}_{0}q(t)b(t)\phi_{1,qb}(t)u_{0}(t)dt \nonumber \\
& < \mu \lambda_{1,qb}\int^{1}_{0}q(t)b(t)\phi_{1,qb}(t)u_{0}(t)dt, \nonumber
\end{align}
a contradiction. Hence our claim holds, that is, (\ref{e20}) has no solution with norm r. Thus, the condition (b) of Theorem \ref{thm2.2} is satisfied.

\ \ \ \ Now, we consider the set
\[D = \{u \in K ; \, \, r\leq \|u\|\leq R\}.\]
Then clearly, $T:D \to K$ is compact and continuous. We shall prove the condition (a) of Theorem 2.2. To prove this, it is enough to show that for any $\bar{R} \geq R$, the problem $u = \mu Tu, \, \mu > 1$ has no solution of norm $\|\bar{R}\|.$ If this is not true, then there exists a sequence $\{R_{n}\}^{\infty}_{n=1}, \, R_{n} \to \infty$ on $n \to \infty, \, R_{n}\geq \bar{R}$ and a sequence $\{\mu_{n}\}^{\infty}_{n=1}$ of reals with $\mu_{n} > 1$ and a sequence of function $\{u_{n}\}^{\infty}_{n=1}$ with $\|u_{n}\| = R_{n}$ such that $u_{n} = \mu_{n} Tu_{n}$ holds, that is
\begin{equation}\label{e23}
\begin{cases}
-u''_{n}(t) = \mu_{n}\lambda q(t)f(t,u_{n}(t)), \, 0<t_{n}<1,\\
u_{n}(0) = 0 = u_{n}(1).
\end{cases}
\end{equation}
Multiplying the equation $-u''_{n}(t) = \mu_{n}\lambda q(t)f(t,u_{n}(t))$ by $\phi_{1,qb}(t),$ and integrating from $0$ to $1$, we obtain
\begin{align}
-\int^{1}_{0}u''_{n}(t)\phi _{1,qb}(t)dt = & \mu_{n}\lambda \int^{1}_{0}q(t)\phi_{1,qb}(t)f(t,u_{n}(t))dt \nonumber \\
& > \mu_{n}c\lambda \int^{1}_{0}q(t)b(t)\phi_{1,qb}(t)u_{n}(t)dt, \nonumber
\end{align}
that is,
\begin{align}
\mu_{n}c\lambda \int^{1}_{0}q(t)b(t)\phi_{1,qb}(t)u_{n}(t)dt &  < -\int^{1}_{0}u''_{n}(t)\phi _{1,qb}(t)dt \nonumber \\
= & -\int^{1}_{0}u_{n}(t)\phi''_{1,qb}(t)dt \nonumber \\
= & \lambda_{1,qb}\int^{1}_{0}q(t)b(t)\phi_{1,qb}(t)u_{n}(t)dt. \label{e24}
\end{align}
Since $\lambda > \frac{\lambda_{1,qb}}{c}$ and $\mu_{n}>1,$ then (\ref{e24}) yields a contradiction. Hence the condition (a) of Theorem \ref{thm2.2} is satisfied. The proof of condition (c) of Theorem \ref{thm2.2} follows from Lemma \ref{Lemma 2.4}. By Theorem \ref{thm2.2}, the BVP(\ref{e4}) has a positive solution in D. The theorem is proved.\\

\noindent    \textbf{\underline{Proof of Theorem \ref{thm4.2}}:}
We shall use Theorem \ref{thm2.1} to prove the theorem. Let $r \in (0,\delta).$ Let $u(t)$ be a solution of $u = \mu Tu$ with $\mu > 1.$ We claim that $\|u\| \neq r.$ If this is not true, there exists a solution $u_{0}(t)$ of $u(t) = \mu Tu(t), \, \, \mu > 1,$ and $u_{0}(t)$ satisfies the property $\|u_{0}\| = r.$ Then $u_{0}(t)$ is a solution of
\begin{equation}\label{e25}
u''_{0}(t) + \lambda \mu q(t)f(t,u_{0}(t)) = 0,\, \, 0<t<1, \, \mu>1
\end{equation}
together with the boundary condition
\[u_{0}(0) = u_{0}(1) = 0.\]
Multiplying both sides of Eq.(\ref{e25}) by $\phi_{1,qb}(t)$ and integrating from $0$     to $1$, we obtain, using (H4) and $\lambda > \frac{\lambda_{1,qb}}{c}$, that
\begin{align}
\mu \lambda_{1,qb}\int^{1}_{0}q(t)b(t)\phi_{1,qb}(t)u_{0}(t)dt = & -\int^{1}_{0}u''_{0}(t)\phi_{1,qb}(t)dt  \nonumber \\
= & -\int^{1}_{0}u_{0}(t)\phi''_{1,qb}(t)dt \nonumber \\
= &  \lambda_{1,qb}\int^{1}_{0}q(t)b(t)\phi_{1,qb}(t)u_{0}(t)dt, \nonumber
\end{align}
which is a contradiction. Hence our claim holds.   Thus, if we consider the set
\[D = \{u \in K ; r\leq \|u\|\leq R\},\]
then $T: D \to K$ is compact and continuous. Further, for the above choice of $r,$ the condition (b) of Theorem \ref{thm2.1} is satisfied.

\ \ \ \ Now, we prove the condition (a) of Theorem \ref{thm2.1}. Let $u(t) \in D$ be a solution of $u = \mu Tu$ and $\mu<1.$ We shall show that $\|u\| \neq R.$ For this, it is enough to show that the problem $u= \mu Tu, \, \mu < 1$ has no solution of norm $\bar{R}$ for any $\bar{R} \geq R.$ If possible, suppose that there exists a solution $u_{1}(t)$ of $u = \mu Tu, \, \mu<1$ such that $\|u_{1}\| = R_{0}, \, \, R_{0} \geq R.$ Since $u_{1}(t)$ is a solution of
\begin{equation}\label{e26}
u''_{1}(t) + \lambda\mu q(t)f(t,u_{1}(t)) = 0,\, \, 0<\mu<1
\end{equation}
with
\[u_{1}(0) = 0 = u_{1}(1),\]
then multiplying both sides of (\ref{e26}) by $\phi_{1,qb}(t),$ integrating from $0$ to $1,$ and using $\lambda < \lambda_{1,qb},$ we have
\begin{align}
\lambda_{1,qb}\int^{1}_{0}q(t)b(t)\phi_{1,qb}(t)u_{1}(t)dt = & \mu\lambda\int^{1}_{0}q(t)\phi_{1,qb}(t)f(t,u_{1}(t))dt \nonumber  \\
& \leq \mu\lambda\int^{1}_{0}q(t)\phi_{1,qb}(t)b(t)u_{1}(t)dt\nonumber \\
& < \lambda\int^{1}_{0}q(t)\phi_{1,qb}(t)b(t)u_{1}(t)dt, \nonumber
\end{align}
a contradiction. Hence our claim holds, which proves the condition (a) of Theorem \ref{thm2.1}. \\
The proof of the condition (c) of Theorem \ref{thm2.1} is similar to the proof of Lemma \ref{Lemma 2.4}. By Theorem \ref{thm2.1}, the BVP (\ref{e4}) has atleast one positive solution $u(t)$.  This completes the proof of the theorem.\\

\end{document}